\documentclass[a4paper]{article}
\usepackage{amsmath, amssymb, amsfonts, amscd, amsthm}

\newcommand{\Z}{\mathbb{Z}}

\renewcommand{\o}{\omega}
\renewcommand{\l}{\Lambda}
\newcommand{\p}{\varphi}
\renewcommand{\t}{\theta}
\newcommand{\<}{\langle}
\renewcommand{\>}{\rangle}
\renewcommand{\top}{\mathfrak{Top}}

\newtheoremstyle{pedro}{}{}{\itshape}{}{\sc}{~--}{ }{}

\theoremstyle{pedro}
\newtheorem{lem}{Lemma}[section]

\newtheorem{thm}[lem]{Theorem}

\newtheorem{prop}[lem]{Proposition}

\newtheorem{coro}[lem]{Corollary}

\theoremstyle{remark}

\theoremstyle{definition}

\usepackage{titlesec}

\titleformat{\section}{\bf\center}{\S \arabic{section}.}{1em}{}
\titleformat{\subsection}{\bf}{}{0pt}{}

\title{The representation ring of a simply connected Lie group as a $\lambda$-ring}
\begin{document}

\maketitle

\begin{abstract}Adams and Conway have stated without proof a result which says, roughly speaking, that the representation ring $R(G)$ of a compact, connected Lie group $G$ is generated as a $\lambda$-ring by elements in 1-to-1 correspondance with the branches of the Dynkin diagram. In this note we present an elementary proof of this.
\end{abstract}

\section{Introduction}

Let $G$ be a $1$-connected, compact Lie group. A celebrated theorem of Weyl's describes the representation ring $R(G)$ as a polynomial ring $$R(G)=[V_1,\ldots,V_r],$$ where the "fundamental representations" $V_i$ are in 1-to-1 correspondance with the nodes of the Dynkin diagram. Therefore, in order to understand the representations of $G$, one basically needs to understand $r$ of them.

In his 1985 article \cite{adamse8}, J F. Adams presents a method, which he ascribes to J. H. Conway, to cut down on this number. The idea, simply enough, is to exploit the structure of $\lambda$-ring which $R(G)$ possesses since it is possible to take exterior powers of representations (see \cite{atiyahtall} for this -- we shall not need anything from the general theory of $\lambda$-rings, though).

The key result announced by Adams is the following. Suppose that the Dynkin diagram of $G$ has an arm of length $k$; by this we mean that there are nodes $v_1,v_2,\ldots,v_k$ such that $v_i$ is connected to $v_{i+1}$ by a single bond for $1\le i < k$, and such that no other edge runs to any of $v_1, v_2, \ldots, v_{k-1}$. Then one has:

\begin{thm}\label{thm:main} Suppose that the Dykin diagram of the $1$-connected, compact Lie group $G$ has an arm of length $k$. Let $V_i$ be the representation corresponding to $v_i$. Then for $1\le i \le k$, one has $$\l^i V_1=V_i \oplus \mathrm{(lower~terms)}$$ where the lower terms are irreducible representations $W < V_i$.
\end{thm}

In this statement we order the irreducible representations of $G$ according to their respective maximal weights.

For any $G$-module $M$, we let $\top(M)$ denote the smallest submodule of $M$ having the same maximal weights, with the same multiplicities. In other words, $\top(M)$ is the smallest $G$-submodule of $M$ satisfying $M=\top(M) \oplus \mathrm{(lower~terms)}$ (here we mean that any irreducible factor in the "lower terms" is dominated by at least one irreducible factor in $\top(M)$). With this terminology, theorem \ref{thm:main} can be reformulated thus:
$$V_i=\top(\l^iV_1).$$

As another illustration of the $\top$ notation, which is due to Adams, we point out that the traditional proof of Weyl's theorem (eg as in \cite{adams2}) shows that any irreducible representation of $G$ is of the form $\top(V_1^{\otimes n_1}\otimes V_2^{\otimes n_2}\otimes\cdots\otimes V_r^{\otimes n_r})$.

With this result in hand, one can show easily that $\l^iV_1$ may replace $V_i$ in the statement of Weyl's theorem (this is more or less trivial, but requires a certain amount of notation -- we shall come back to this in due time). Therefore, if we see $R(G)$ as a $\lambda$-ring, the whole arm of the Dynkin diagram is accounted for by the single generator $V_1$.

Adams (unlike Conway) had a proof of theorem \ref{thm:main} ready, but wanted to check if anyone else had previously obtained the result before publishing it. His intention seemed to include the proof in a book he had in preparation on exceptional Lie groups. This was sadly prevented by Adams's untimely death. The book has appeared \cite{adams}, based on his lecture notes, but does not include this proof.

In this quick note we present an elementary proof of the main theorem, and give a few applications.

\section{Notations}

$G$ is a $1$-connected, compact Lie group, with maximal torus $T$ and Weyl group $W$. We fix a $W$-invariant inner product on $L(T)$, the real Lie algebra tangent to $T$, denoted by $\<-,-\>$.

The roots of $G$ are non-zero weights of the adjoint representation. They come in $m$ pairs $\pm\t_1,\pm \t_2,\ldots,\pm\t_m$, and we assume that the signs are arranged so that the "fundamental Weyl chamber" $$FWC=\{v\in L(T): \t_i(v)>0 \}$$ is nonempty. We call $\t_1,\ldots,\t_m$ the positive roots; we assume that the indices are so chosen that the simple roots are $\t_1,\t_2,\ldots,\t_r$ (where $r=\dim T$).

 Dually, the fundamental dual Weyl chamber is $$FDWC=\{\lambda\in L(T)^*: \<\lambda,\t_i\> > 0 \}.$$
To each simple root $\t_i$ corresponds the fundamental weight $\o_i$ characterized by
$$\frac{2\<\o_i,\t_j\>}{\<\t_j,\t_j\>}=\delta_{ij}$$ for any simple root $\t_j$. 
 
We partially order weights by requiring $\lambda_1 < \lambda_2$ if and only if $\lambda_1(v) < \lambda_2(v)$ for all $v\in FWC$.

For any root $\t_i$, there is an element $\p_i\in W$ which induces on $L(T)$ the orthogonal reflexion in $\ker \t_i$. For convenience, $\p_0$ will denote the identity. Very explicitly, for $\lambda\in L(T)$ or $\lambda\in L(T)^*$:

$$\p_i(\lambda)=\lambda - \frac{2\<\lambda,\t_i\>}{\<\t_i,\t_i\>}\t_i.$$

Finally, we assume, of course, that the Dynkin diagram of $G$ has an arm of length $k$, as explained in the introduction. For simplicity, we arrange the notations so that $\t_1,\t_2,\ldots,\t_k$ are the simple roots corresponding to the nodes $v_1,v_2,\ldots,v_k$. We write $V_i$ for the representation associated with $\o_i$ (or with $v_i$, if you want). We will also write $V$ for $V_1$ and $\o$ for $\o_1$.

It may be useful to recall that two simple roots $\t_i$ and $\t_j$ correspond to nodes with a simple bond between them precisely when $\<\t_i,\t_i\>=\<\t_j,\t_j\>$ and $$\frac{-2\<\t_i,\t_j\>}{\<\t_i,\t_i\>}=1.$$
For example, we shall use repeatedly $\p_i\t_{i+1}=\p_{i+1}\t_i=\t_i + \t_{i+1}$ for $i<k$.

\section{Main proof}

Our first objective is to show that the weight $\o_i$ occurs in $\l^i V$.

\begin{lem}\label{lem:labase} For $1\le i < k$, we have
$$\o_{i+1}  =  \o + \p_1\o + \p_2\p_1\o + \ldots + \p_i\p_{i-1}\cdots\p_1\o.$$
\end{lem}

\begin{proof} Let $\lambda$ denote the right hand side. We shall compute $\<\lambda,\t_j\>$.

If $j>i+1$ then we have
$$\<\p_s\cdots\p_1\o,\t_j\>=\<\o,\p_1\cdots\p_s\t_j\>=\<\o,\t_j\>=0$$ for $s\le i$ and it follows that $\<\lambda,\t_j\>=0$.

If now $j<i+1$, the equations above still hold when $s<j-1$. The two terms in $\<\lambda,\t_j\>$ obtained for $s=j-1$ and $s=j$ cancel each other, for $\p_j\t_j=-\t_j$. Finally, for $s\ge j+1$ we obtain
$$\<\o,\p_1\cdots\p_s\t_j\>=\<\o,\p_1\cdots\p_{j+1}\t_j\>=\<\o,\t_{j+1}\>=0$$ where we use $\p_j\p_{j+1}(\t_j)=\p_j(\t_j + \t_{j+1})=-\t_j + \t_j + \t_{j+1}=\t_{j+1}$. So $\<\lambda,\t_j\>=0$ in this case also.

It remains to compute $\<\lambda,\t_{i+1}\>$. Arguing as above, we see that all terms vanish except
$$\<\o,\p_1\cdots\p_i\t_{i+1}\>=\<\o,\t_1 + \t_2 + \cdots + \t_{i+1}\>=\<\o,\t_1\>$$ and we obtain $\<\lambda,\t_{i+1}\>=\<\o,\t_1\>=\frac{1}{2}\<\t_1,\t_1\>=\frac{1}{2}\<\t_{i+1},\t_{i+1}\>$. It follows that $\lambda=\o_{i+1}$.
\end{proof}

As a corollary, the weights $\p_i\p_{i-1}\cdots\p_1\o=\o_{i+1} - \o_i$, for various values of $i$, are pairwise distinct. It follows that, if $v\in V$ is a non-zero eigenvector under the action of $T$ with weight $\o$, then $v\wedge\p_1 v\wedge\ldots\wedge\p_i\cdots\p_1v\in \l^{i+1}V$ is a non-zero eigenvector with weight $\o_{i+1}$.

Let us prove another technical lemma before proceeding.

\begin{lem}
Let $i_1, i_2, \ldots i_t$ be integers, where $1\le t <k$. Then
$$\p_{i_1}\p_{i_2}\cdots\p_{i_t}\o=\p_s\p_{s-1}\cdots\p_1\o$$ for some $s\le t$.

More precisely, for $s<k-1$, the weight $\p_s\p_{s-1}\cdots\p_1\o$ is fixed by any $\p_u$ unless $u=s$ or $u=s+1$.
\end{lem}

\begin{proof} Clearly the first equality follows by repeated use of the more precise statement. Let $\lambda=\p_s\p_{s-1}\cdots\p_1\o=\o_{s+1} - \o_s$. Then
$$\<\p_u\lambda,\t_j\>=\<\lambda,\p_u\t_j\>=\<\lambda,\t_j\> + c_j\<\lambda,\t_u\>$$ for some constant $c_j$. If $u$ is neither equal to $s$ nor $s+1$, then $\<\o_{s+1} - \o_s,\t_u\>=0$ and we conclude that $\<\p_u\lambda,\t_j\>=\<\lambda,\t_j\>$ for all $j$. 
\end{proof}

It follows in particular that $\p_i\o=\o$ if $i\ge 2$. An immediate induction then gives $$\p_i\p_{i-1}\cdots\p_1\o=\o - (\t_1 + \t_2 + \ldots + \t_i)$$ which will be handy in the sequel.

We are now ready to prove that $\o_i$ is maximal in $\l^iV$. Let us start by ordering the weights in the $W$-orbit of $\o$.

\begin{lem}\label{lem:order} Let $v\in FWC$. Then 
$$\o(v)> \p_1\o(v) > \p_2\p_1\o(v) >\cdots >\p_{k-1}\p_{k-2}\cdots\p_1\o(v)$$ and $\p_{k-1}\p_{k-2}\cdots\p_1\o(v) > \p\o(v)$ for any other value of $\p\o(v)$ with $\p\in W$.
\end{lem}

\begin{proof} The first inequalities are trivial, and follow from the last expression given for $\p_i\p_{i-1}\cdots\p_1\o$.
% they follow from
%$$\p_u\lambda=\lambda - \frac{2\<\lambda,\t_u\>}{\<\t_u,\t_u\>}\t_u$$ for any weight $\lambda$, while $\t_u(v)> 0$ for $v\in FWC$ and $\<\o_s - \o_{s-1},\t_s\> > 0$.

It is well-known (\cite{adams2}, 6.26) that $\p\o(v)<\o(v)$ as soon as $\p\o\ne\o$. Now, let $\p\in W$ be such that $\p\o(v)$ reaches its second highest value. Then $\p\o\ne\o$ so $\p\o$ is not in $Cl FDWC$, and therefore there is an $u$ such that $\<\p\o,\t_u\> < 0$. Compute then 
$$\p_u\p\o(v)=\p\o(v) - \frac{2\<\p\o,\t_u\>}{\<\t_u,\t_u\>}\t_u(v) > \p\o(v).$$ By definition of $\p$, this implies $\p_u\p\o(v)=\o(v)$, so that $\p_u\p\o=\o$ and $\p\o=\p_u\o$. By the last lemma, $\p_u\o$ can only be $\p_1\o$, so we conclude that $\p_1\o(v)$ is the second highest value of $\p\o(v)$ for $\p\in W$, and that it is only reached when $\p_1\p$ fixes $\o$.

We would prove similarly that $\p_2\p_1\o(v)$ is the third highest value of $\p\o(v)$, which is only reached when $\p_2\p_1\p$ fixes $\o$. The result follows by induction.
\end{proof}

In order to deal with the weights in $V$ which are not in the $W$-orbit of $\o$, we shall need\footnote{such weights do exist, cf the $26$-dimensional representation of $F_4$ which has two zero weights.} the following.

\begin{lem}\label{lem:adams} Let $M$ be an irreducible $G$-module with maximal weight $\lambda$.
\begin{enumerate}
\item Let $m\in M$ be an eigenvector for the weight $\alpha$. For each positive root $\t_i$, let $X_{\t_i}\in L(G)$ be an eigenvector for $\t_i$. Then $\alpha=\lambda$ if and only if $X_{\t_i}m=0$ for each $i$.
\item All the weights occuring in $M$ are of the form $$\lambda - \sum n_i\t_i,$$ where $n_i$ is a nonnegative integer. 
\end{enumerate}
\end{lem}

\begin{proof} Half of (1) is trivial: if $m$ is an eigenvector for $\lambda$, then $X_{\t_i}m$ is an eigenvector for $\lambda + \t_i > \lambda$, so it must be zero. We prove the converse and (2) at the same time.

Suppose then that $X_{\t_i}m=0$ for each positive root $\t_i$. Since $M$ is irreducible, it is generated as an $L(G)$-module by $m$ alone. Therefore $M$ is spanned by elements of the form $$X_sX_{s-1}\cdots X_1m$$ with $X_j\in L(G)$. Splitting $L(G)$ according to the adjoint action of the maximal torus $T$, we may assume that each $X_j$ is an eigenvector associated to the root $\lambda_j$, or to the zero weight.  A vector as above is then an eigenvector in $M$ with weight $\alpha + \lambda_1 + \ldots \lambda_s$. We claim that $M$ is in fact generated by elements as above for which each $\lambda_i$ is a negative root (or zero). Granting this, $\alpha$ is clearly maximal, so $\alpha=\lambda$, and the other weights in $M$ are as announced in (2).

To prove the claim, we show by induction on $s$ how to replace the generating vectors. The case $s=1$ is our assumption, while the induction proceeds easily using $X_iX_j=X_jX_i + [X_i,X_j]$.
\end{proof}

Theorem \ref{thm:main} now follows from this:

\begin{prop} $\o_{i+1}$ is the highest weight is $\l^{i+1}V$ (equivalently, it is the only maximal weight).
\end{prop}

\begin{proof} By induction on $i$ (the case $i=0$ being given). So let $$\lambda=\lambda_1 + \lambda_2 + \ldots + \lambda_{i+1}$$ be a maximal weight occuring in $\l^{i+1}V$, where each $\lambda_j$ is a weight in $V$.

{\em First case.} Suppose to start with that for each $\lambda_j$ and each $v\in FWC$ we have $\lambda_j(v)\ge \p_i\p_{i-1}\cdots\p_1\o(v)$. By the preceding lemma, we have
$$\lambda_j=\o - \sum_s n_s\t_s$$ while $\p_i\p_{i-1}\cdots\p_1\o=\o - \t_1 - \t_2 - \ldots -\t_i$. It follows that $$\sum_s n_s\t_s(v) \le \t_1(v) + \ldots \t_i(v)$$ for all $v\in FWC$, or even for $v\in Cl FWC$. Applying this with $v=\o_s$ (and identifying $L(G)$ with its dual) we conclude that each $n_s$ is $0$ or $1$ when $1\le s \le i$ (and zero otherwise).

So if $\lambda_j\ne\o$, let $t$ denote the largest value of $s$ such that $n_s=1$, and write $$\lambda_j=\p_t\p_{t-1}\cdots\p_1\o + \sum_{s<t} m_s \t_s,$$ with $m_s=1-n_s$. If there exists an eigenvector in $V$ for this weight, then there is also one for the weight $$\alpha=\o + \p_1\p_2\cdots\p_t\left( \sum_{s<t} m_s\t_s \right).$$ Here $\p_1\p_2\cdots\p_t \t_s=\p_1\cdots\p_{s+1}\t_s=\t_{s+1}$ (as already observed during the proof of \ref{lem:labase}). It follows that $\alpha>\o$ if any coefficient $m_s$ is nonzero. As this is impossible, we conclude that
$\lambda_j=\p_t\p_{t-1}\cdots\p_1\o$.

It follows that $t$ must be different for each $j$, and that $\lambda=\o_{i+1}$.

{\em Second case.} Suppose now that there is a $v\in FWC$ and a $j$ such that $\lambda_j(v)<\p_i\p_{i-1}\cdots\p_1\o(v)$, say $j=i+1$. We shall derive a contradiction (recall that $\lambda$ is assumed to be maximal).

Now if $\lambda_{i+1}$ is in the $W$-orbit of $\o$, we draw from \ref{lem:order} that $\lambda_j(v)<\p_i\p_{i-1}\cdots\p_1\o(v)$ actually holds for {\em all} $v\in FWC$. Then by the induction hypothesis, we have
$$\lambda_1(v) + \lambda_2(v) + \ldots \lambda_i(v) \le \o(v) + \p_1\o(v) + \ldots + \p_{i-1}\p_{i-2}\cdots\p_1\o(v),$$ so $\lambda(v)<\o_{i+1}(v)$ and we are done.

If on the other hand $\lambda_{i+1}$ is not in the $W$-orbit of $\o$, we shall proceed quite differently. Let $x\in V$ be an eigenvector associated to $\o$, and let $y\in V$ be an eigenvector for $\lambda_{i+1}$. Put
$$m=x\wedge \p_1x \wedge\ldots \wedge \p_{i-1}\p_{i-2}\cdots\p_1x\wedge y\in \l^{i+1}V.$$
Then $m$ is non-zero, and is an eigenvector for a weight which is $\ge \lambda$ (using the induction hypothesis), so by maximality of $\lambda$, we see that $m$ is an eigenvector for $\lambda$. Using lemma \ref{lem:adams}, we conclude that there exists a positive root $\t_j$ and an eigenvector $X=X_{\t_j}\in L(G)$ such that $Xy\ne 0$. However by the same lemma, $Xx=0$ and $Xm=0$ ($m$ being an eigenvector for the maximal weight of some irreducible submodule of $\l^{i+1}$). Write
$$\begin{array}{rcl}
0= Xm & = & X(x\wedge \p_1x \wedge\ldots \wedge \p_{i-1}\p_{i-2}\cdots\p_1x)\wedge y \\ &   & \pm x\wedge \p_1x \wedge\ldots \wedge \p_{i-1}\p_{i-2}\cdots\p_1x\wedge Xy
\end{array}$$
The first term here is $0$ by induction, so we conclude that $Xy$ is in the linear span of $x, \p_1x, \ldots, \p_{i-1}\p_{i-2}\cdots\p_1x$. Since $Xy$ is a (nonzero) eigenvector for the action of $T$, it follows that $Xy$ is in fact (a scalar multiple of) one of these. In conclusion, $\lambda_{i+1} + \t_j=\p_s\p_{s-1}\cdots\p_1\o$ for some $s$ with $1\le s \le i-1$, or equivalently $\lambda_{i+1}=(\o_{s+1} - \o_s) - \t_j$.

Now, the contradiction will come from the fact that, by maximality of $\lambda$, we must have $\lambda_{i+1}\in Cl FDWC$, so that $\p_j\lambda_{i+1} \le \lambda_{i+1}$. This amounts to
$$ \frac{\<\o_{s+1} -\o_s,\t_j\>}{\<\t_j,\t_j\>} \ge 1.$$ For this to happen, we would need a positive root $\t_j$ which, when written as a sum of positive simple roots $\t_j=\sum n_t\t_t$, has $n_{s+1}-n_s \ge 2$. This never happens. You may either inspect the existing root systems, or fill in the details of the following sketch: let $S$ be the set of positive roots such that $|n_{s+1}-n_s| \le 1$ for all $1\le s+1 \le k-1$. Then $S$ contains the positive simple roots, and is stable under each relexion $\p_u$ for $u\le r$. Hence $S$ is stable under $W$ and consequently, it is equal to the set of all positive roots.
%% Now, by maximality of $\lambda$, we see that $\lambda_{i+1}\in Cl FDWC$, so %that $\p_j\lambda_{i+1} < \lambda_{i+1}$. This amounts to
%$$ \frac{\<\o_{s+1} -\o_s,\t_j\>}{\<\t_j,\t_j\>} > 1.$$ If $\t_j$ were a simple root, this would imply $\t_j=\t_s$. In general, one can see from the shape of the Dynkin diagram that $\t_j$ either appears in the following list: $\t_s, \t_s + \t_{s+1}, \ldots, \t_s + \t_{s+1} + \ldots + \t_i$, or is of the form $\t_s + \t_{s+1} + \ldots \t_i + \t$, where $\t$ is a sum of positive simple roots $\t_u$ with $u>i$. We rule out this possibility at once: it would imply that $\lambda_{i+1}=(\o_{s+1} - \o_s) - \t_j=\o - (\t_1 + \ldots + \t_i + \t)$ is dominated by $\o_{i+1} - \o_i= \o - (\t_1 + \ldots + \t_i)$, which contradicts the maximality of $\lambda$.
%
%So we must have $\lambda_{i+1}=\o - (\t_1 + \ldots + \t_t)=\o_{t+1} - \o_t$ for some $t\le i$. This is also a contradiction, since there is a $v\in FWC$ such that $\lambda_{i+1}(v)<\o_{i+1}(v) - \o_i(v) \le \o_{t+1}(v) - \o_t(v)$ (using lemma \ref{lem:order}).
\end{proof}

\begin{coro} We may replace $V_i$ by $\l^iV$ in the statement of Weyl's theorem, for $1\le i \le k$. More generally, if $\top (V_i')=V_i$ for each $i$, then $$R(G)=\Z[V_1',V_2',\ldots,V_r']$$ and moreover any irreducible representation of $G$ is of the form $$\top({V'_1}^{\otimes n_1}\otimes {V'_2}^{\otimes n_2}\otimes\cdots\otimes {V'_r}^{\otimes n_r}).$$
\end{coro}

\begin{proof} Run through a proof of Weyl's theorem and replace $V_i$ by $V_i'$.

As this is short enough to write out, we give the details. For each weight $\lambda\in Cl FDWC$, let $S(\lambda)$ denote the elementary symmetric sum $$S(\lambda)=\sum_{\alpha\in W\lambda} e^{2i\pi \alpha}.$$ These form a $\Z$-basis for $R(G)=R(T)^W$ (viewed as the character ring rather than the representation ring).

Now, if $\lambda=n_1\o_1 + \ldots + n_r\o_r$, we clearly have
$$\chi({V'_1}^{\otimes n_1}\otimes\cdots\otimes {V'_r}^{\otimes n_r})=S(\lambda) + \sum_{\lambda_i<\lambda} m_i S(\lambda_i)$$ where $\chi$ means the character (use \cite{adams2}, 6.36). If we suppose by induction that each $S(\lambda_i)$ can be written in the form $\chi(P(V_1',\ldots,V_r'))$ (with $P$ a polynomial), then the same can be said of $S(\lambda)$.

As a result, $V_1', V_2',\ldots, V_r'$ generate $R(G)$ as a ring. They are clearly algebraically independant (for example, consider the transcendance degree of $R(G)$).  

The rest is easy.
\end{proof}

\section{Applications}

All the information on exceptional Lie groups used below can be found in Adams's book\cite{adams}.

\subsection{The representation ring of $F_4$.} This group has the following Dynkin diagram:

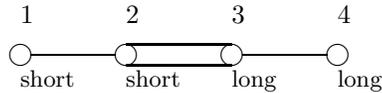
\setlength{\unitlength}{4pt}

\begin{center}
\begin{picture}(45,5)
\multiput(0,0)(10,0){4}{\circle{2}}
\put(0,3){1}
\put(10,3){2}
\put(20,3){3}
\put(30,3){4}
\put(0,-3){\small short}
\put(10,-3){\small short}
\put(20,-3){\small long}
\put(30,-3){\small long}
\put(1,0){\line(1,0){8}}
\put(21,0){\line(1,0){8}}
\put(10,1){\line(1,0){10}}
\put(10,-1){\line(1,0){10}}
\end{picture}
\end{center}

\bigskip
The group $F_4$ may be defined as the group of automorphisms of the (Jordan) algebra $J$ of $3\times 3$ hermitian matrices over the octonions, the multiplication being $A,B\mapsto \frac{1}{2}(AB + BA)$. The real dimension of $J$ is $27$, and $F_4$ acts on the subspace of matrices with trace $0$. Complexifying, this affords a $26$-dimensional representation $U$ of $F_4$. This $U$ turns out to be irreducible and its highest weight corresponds to the node $1$.

The (complexified) adjoint representation $Ad$, on the other hand, corresponds to the node $4$.

It follows from the main theorem in this paper that
$$R(F_4)=\Z[Ad,\l^2Ad, U, \l^2 U].$$

\subsection{The $E$ family.}

Adams originally considered the case of $E_8$ in \cite{adamse8}. This group has the following Dynkin diagram:

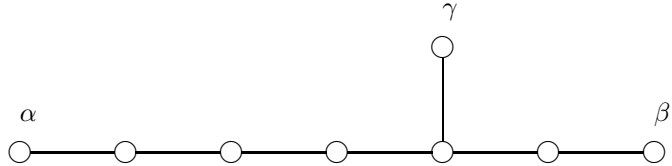
\begin{center}
\begin{picture}(75,15)
\multiput(0,0)(10,0){7}{\circle{2}}
\multiput(1,0)(10,0){6}{\line(1,0){8}}
\put(40,1){\line(0,1){8}}
\put(40,10){\circle{2}}
\put(0,3){$\alpha$}
\put(40,13){$\gamma$}
\put(60,3){$\beta$}
\end{picture}
\end{center}

Adams concludes that
$$R(E_8)=\Z[\alpha,\l^2\alpha,\l^3\alpha,\l^4\alpha,\beta,\l^2\beta,\gamma,\delta]$$ where $\delta$ can be taken to be either $\l^5\alpha$, or $\l^3\beta$, or $\l^2\gamma$. It is easy to see that the representation $\alpha$ is the adjoint representation, and Adams in {\em loc cit} gives some concrete information on $\beta$ and $\gamma$ (eg, the way they sit in $\alpha\otimes\alpha$). I do not know a direct, concrete definition of $\beta$ and $\gamma$.

Similarly, for $E_7$ we have the following diagram:

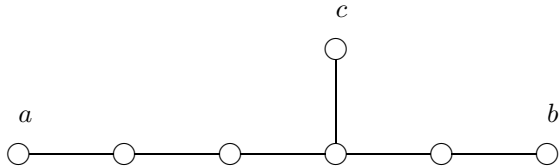
\begin{center}
\begin{picture}(65,15)
\multiput(0,0)(10,0){6}{\circle{2}}
\multiput(1,0)(10,0){5}{\line(1,0){8}}
\put(30,1){\line(0,1){8}}
\put(30,10){\circle{2}}
\put(0,3){$a$}
\put(30,13){$c$}
\put(50,3){$b$}
\end{picture}
\end{center}

So that
$$R(E_7)=\Z[a,\l^2a,\l^3a,b,\l^2b,c,d]$$ where $d$ can be taken to be either $\l^4a$, or $\l^3b$, or $\l^2c$. Here however, the adjoint representation does not correspond to any of $a$, $b$, $c$. We can identify nevertheless the representation $a$ as the complexification of the $56$-dimensional real representation $W$ of $E_7$ which can be used to define $E_7$ as the group of linear automorphisms of $W$ preserving both a quadratic form and a quartic form. 

The situation for $E_6$ is analogous. The Dynkin diagram has three branches, so $R(E_6)$ is generated by $3$ elements as a $\lambda$-ring. I am not aware of any explicit description of the three fundamental representations involved. They are neither given by the adjoint representation, nor by the representation $W$ used to define $E_6$ as a group of maps.

%\begin{picture}(55,15)
%\multiput(0,0)(10,0){5}{\circle{2}}
%\multiput(1,0)(10,0){4}{\line(1,0){8}}
%\put(20,1){\line(0,1){8}}
%\put(20,10){\circle{2}}
%\end{picture}
%

\bibliography{myrefs}

\newcommand{\noopsort}[1]{} \newcommand{\printfirst}[2]{#1}
  \newcommand{\singleletter}[1]{#1} \newcommand{\switchargs}[2]{#2#1}
\begin{thebibliography}{1}

\bibitem{adams2}
{\sc J.~F. Adams}, {\em Lectures on {L}ie groups}, W. A. Benjamin, Inc., New
  York-Amsterdam, 1969.

\bibitem{adamse8}
\leavevmode\vrule height 2pt depth -1.6pt width 23pt, {\em The fundamental
  representations of ${E}_8$}, Contemp. Math, 37 (1985).

\bibitem{adams}
\leavevmode\vrule height 2pt depth -1.6pt width 23pt, {\em Lectures on
  exceptional {L}ie groups}, Chicago Lectures in Mathematics, University of
  Chicago Press, Chicago, IL, 1996.

\bibitem{atiyahtall}
{\sc M.~F. Atiyah and D.~O. Tall}, {\em Group representations, {$\lambda
  $}-rings and the {$J$}-homomorphism}, Topology, 8 (1969), pp.~253--297.

\end{thebibliography}
\bibliographystyle{siam}
\end{document}